\documentclass[11pt]{article}
\usepackage{amsmath, amssymb, graphics, epsfig, enumerate, amstext, amsthm, comment}
\usepackage{amssymb,amsmath}
\usepackage{latexsym,bm,mathrsfs}

\usepackage{graphicx,amssymb,amsmath,latexsym,amsfonts,amsthm}
\usepackage{amsmath,amssymb,amsthm}
\usepackage{latexsym}
\usepackage{amssymb}
\usepackage{stmaryrd}
 \usepackage{float}
\usepackage{psfrag}
\usepackage{xcolor}
\usepackage{enumerate}
\usepackage{manfnt,mathrsfs}
\usepackage{hyperref}
\usepackage{url}

\makeatother
\newtheorem{theorem}{Theorem}[section]

\newtheorem{remark}[theorem]{Remark}

\newtheorem{proposition}[theorem]{Proposition}
\newtheorem{lemma}[theorem]{Lemma}

  \def\tr{{\rm tr}\,}

\def\cP{{\cal P}}

\def\Sym{{\rm Sym}}
\def\Alt{{\rm Alt}}

\newcommand\GAP{\textsf{GAP}}
\newcommand\Cliquer{\textsf{Cliquer}}

\renewcommand\le{\leqslant}
\renewcommand\ge{\geqslant}

\newcommand\supp{{\rm Supp}}
\def\qed{{\hfill\vbox{\hrule width 6 pt \hbox{\vrule height 6 pt width 6 pt}}}}
\begin{document}
\openup .6 \jot
\title{\bf Cliques and independent subgroups of the Birkhoff polytope graph}
\author{Zejun Huang\thanks{School of Mathematical Sicences, Shenzhen University, Shenzhen,  518060, P.R. China. (zejunhuang@szu.edu.cn) }, Chi-Kwong Li\thanks{Department of
Mathematics, William \& Mary, Williamsburg, VA 23187, USA.
(Email: ckli@math.wm.edu)},  Eric Swartz\thanks{Department of Mathematics, William \& Mary, Williamsburg, VA 23187, USA. (Email: easwartz@wm.edu),},
Nung-Sing Sze\thanks{Department of Applied Mathematics, The Hong Kong Polytechnic
University, Hung Hom, Hong Kong. (raymond.sze@polyu.edu.hk)} }
\date{}
\maketitle
\begin{abstract}
The Birkhoff polytope  $\Omega_n$ is the polytope of doubly stochastic matrices of order $n$. The Birkhoff polytope graph $G(\Omega_n)$ is the skeleton of $\Omega_n$; it is the Cayley graph whose vertex set consists of the elements of the symmetric group $\Sym(n)$ of degree $n$, where two permutations are adjacent if one equals the product of the other with a cycle.   We study the combinatorial structure of this graph, focusing on its maximal and maximum cliques and on its independent subgroups (subgroups of $\Sym(n)$ whose elements are pairwise nonadjacent in the graph).   We obtain maximal subgroups of $G(\Omega_n)$ and establish both a lower bound and an upper bound for its clique number.  Especially, we prove that if $K$ is a subset of $\Sym(n)$ consisting of 3-cycle permutations such that $\delta_1^{-1}\delta_2$ is a single cycle for all $\delta_1,\delta_2\in K$, then the maximum size of $K$ is $\lfloor (n-1)^2/4\rfloor$, which can be viewed as an Erd\H{o}s-Ko-Rado-type theorem for $\Sym(n)$.
\end{abstract}


\section{Introduction}

Birkhoff's famous theorem \cite{GB} asserts that the convex hull of the set
$\Sym(n)$ of $n\times n$ permutation matrices is the  set $\Omega_n$
of doubly stochastic matrices, i.e., $n\times n$ nonnegative matrices
with each row sum and column sum equal to 1. The set $\Omega_n$ is sometimes referred
to as  the Birkhoff polytope or  the  assignment polytope. Many interesting
problems and results were raised from theory or applications related to the
Birkhoff polytope. For instance, the
van der Waerden conjecture states that
${\rm per}(A)\ge n!/n^n$ for all $A\in
\Omega_n$, with equality when $A$ is the matrix with every entry equal to $1/n$, where
${\rm per}(A)$ is the permanent of $A$. This beautiful conjecture was posed in 1926 \cite{BW}
and had been a long-standing open problem. It was proved independently by Egorychev
\cite{GE} and Falikman \cite{DF} in 1981.

The Birkhoff polytope graph  $G(\Omega_n)$  is the skeleton of $\Omega_n$, whose vertex
set is $\Sym(n)$ and  two vertices in $G(\Omega_n)$ are adjacent if  their convex hull
is an edge of $\Omega_n$.  The adjacency relation admits some alternative
descriptions using graph theory and abstract algebra theory as follows.

For any $A=(a_{ij})\in \Omega_n$, we can define a bipartite
graph $G_B(A)$  with two vertex parts $\{1,2,\ldots,n\}$ such $(i,j)$ is an edge if and
only if $a_{ij}>0$. Then two vertices $A_1$ and $A_2$ in $G(\Omega_n)$
are adjacent if and only if $G_B(A_1)\cup G_B(A_2)$  contains exactly one cycle
\cite{BR2}.

One can regard $\Sym(n)$ as the symmetric group of permutations on
$[n] = \{1, \dots, n\}$. Then two permutations $\sigma_1$ and $\sigma_2$
are adjacent if $\tau = \sigma_1^{-1}\sigma_2$  is a cycle $(i_1, \dots, i_k)$ of length $k > 1$, i.e., $\tau$ is the permutation sending
$i_j$ to $i_{j+1}$  for $j=1, \dots k-1$, sending $i_{k}$ to $i_1$,
and fixing all other numbers in $[n]$.  In this paper, we will use the convention of multiplying permutations right-to-left; that is, if $\sigma_1$ and $\sigma_2$ are permutations, the permutation $\sigma_1 \sigma_2$ represents the permutation obtained by first applying $\sigma_2$, then applying $\sigma_1$.
Denote by $\mathcal{C}_n$ the set of cycles
in $\Sym(n)$. Then $G(\Omega_n)$ is the  Cayley
graph  ${\rm Cay}(\Sym(n),\mathcal{C}_n)$   with vertex set $\Sym(n)$ and edge set
$\{(\delta,\delta\tau): \delta\in \Sym(n),\tau\in \mathcal{C}_n\}$; see \cite{RAB}.

The Birkhoff polytope graph is  vertex transitive, Hamiltonian, and regular
with diameter 2; see \cite{BR1, BR2, RAB}. Moreover, each vertex in $G(\Omega_n)$ has
degree $\sum_{k=2}^n {n \choose k} (k-1)!$ (\cite{BR2}).
Kane, Lovett, and Rao \cite{KLR} studied the independence number $\alpha(n)$ of $G(\Omega_n)$ -- defined to be the size of the largest set of vertices such that no two are adjacent -- which
has applications in coding theory. They
obtained an upper bound $\alpha(n)\le 2(n!/2^{n/2})$. They also presented a large
independent set with cardinality $n!/4^{n}$ when $n=2^m$. Coregliano and Jeronimo
\cite{CJ} improved the upper bound  to be $O(n!/1.97^n)$. By constructing large
independent sets, they also presented a lower bound
$\alpha(n)\ge \Pi_{i=1}^{\lfloor log_2(n)\rfloor}\lfloor n/2^i\rfloor!$.
One may refer to \cite{RAB} for more results on the Birkhoff polytope graph.

A clique of a graph $G$ is a subset $S$ of its vertex set such that the induced subgraph on $S$ is a complete subgraph of $G$.
The clique number of a graph $G$, denoted $\omega(G)$, is the maximum size of a clique in $G$. We denote by $\omega(n)$ the clique number of $G(\Omega_n)$.
For $k\le n$, let $\mathcal{C}_n(k)$ be the  set of all cycles with length $k$ in  $\Sym(n)$.  If a clique $K$ of $G(\Omega_n)$ is a subset of $\mathcal{C}_n(k)$, we call $K$ a {\it $k$-cycle clique}. If $K$ has the largest size among all $k$-cycle cliques, it is called a {\it maximum $k$-cycle clique}. If an independent set of $G(\Omega_n)$ is a subgroup of $\Sym(n)$, we call it an
{\it independent subgroup}.

To determine the exact clique number of $G(\Omega_n)$ is   difficult.  In this article, we study the bounds on $\omega(n)$. We characterize the maximum  $k$-cycle cliques for $k=2,3$,  present some maximal independent subgroups of  $G(\Omega_n)$ and then use them to obtain bounds on $\omega(n)$.

In Section \ref{sect:clique}, we  present a lower bound and an upper bound on $\omega(n)$.
Especially,  we characterize the maximum $k$-cycle cliques of $G(\Omega_n)$ for $k=2,3$ by applying results from the extremal set theory and the extremal graph theory.  In Sections \ref{sect:subgroups}, we  present some maximal independent subgroups of  $G(\Omega_n)$. Some future research directions are
mentioned in Section \ref{sect:related}, and we include in Appendix \ref{appendix} various computational results.

\section{The clique number of $G(\Omega_n)$}
\label{sect:clique}
 In this section, we  study the clique number of $G(\Omega_n)$.
 In particular, we present a  lower bound  and an upper bound  on $\omega(n)$. We also characterize the maximum $k$-cycle cliques, i.e., cliques
 consisting of $k$-cycles, of $G(\Omega_n)$ for $k=2,3$.

Notice that if $X$ is a  maximum clique of $G(\Omega_n)$, then the set $\sigma^{-1} X\equiv\{\sigma^{-1}\tau, \tau\in X\}$ is also a maximum clique of $G(\Omega_n)$. So to find the clique number $\omega(n)$, we may always assume the identity permutation is in the maximum clique. Thus the problem is reduced to finding the largest clique in which each permutation is a single cycle.

Denote by $I_n \in \Sym(n)$   the identity permutation.
It is clear that $G(\Omega_n)$ itself a complete graph when $n < 4$
and thus $\omega(1) = 1$, $\omega(2) = 2$ and $\omega(3) = 6$.
For $n=4$, one knows that the set $\Sym(4)$ can be partitioned into 6 cosets of
$H = \{I_4, (12)(34),(13)(24),(14)(23)\}$. Suppose $R$ is a maximum clique in
$G(\Omega_4)$. Then $R$ contains at most one permutation  from each of these cosets.
Hence $|R|\le 6$. On the other hand, notice that the degree of each vertex in
$G(\Omega_4)$ is 20. Suppose $R$ is a set which contains exactly one permutation from
each  of these cosets. Then $R$ is a clique of  $G(\Omega_4)$. Therefore, we have
$\omega(4)=6.$ By a computer search using \GAP \cite{GAP}, we have $\omega(5)=13$, $\omega(6)=18$, $\omega(7)=23$; see Appendix \ref{app:max567} for maximum cliques in these cases.
So, we have the following.
\begin{center}
\begin{tabular}{|c|c| c| c| c| c| c| c|} 
\hline 
$n$ & 1 & 2 & 3 & 4 & 5 & 6 & 7\cr
\hline
$\omega(n)$ & 1 & 2 & 6 & 6 &13 & 18 & 23\cr\hline
\end{tabular}
\end{center}

For $n \ge 8$, it is not easy to determine $\omega(n)$ by computer search.
We obtain bounds for $\omega(n)$ by studying $k$-cycle cliques, i.e., cliques consisting
of cycles of length $k$ with $k > 1$.
For $\sigma,\tau\in \mathcal{C}_n$, we denote by $\sigma\cap \tau$ the set of
common numbers in $\sigma$ and $\tau$. The following observation is useful in our discussion.

\begin{proposition}\label{p46}
Let $\sigma,\tau$ be two cycles  in $\Sym(n)$ of length at least 3 and
$\sigma\cap \tau=\{1,2\}$.  Then $\sigma \tau$ is a cycle if and only if
$(\sigma, \tau)$ has the form $((1,2,\ldots),(2,1,\ldots))$ or
$((2,1,\ldots),(1,2,\ldots))$.
\end{proposition}

\noindent\it Proof. \rm
 Suppose $\sigma=(1,J_1,2,L_1),\tau=(1,J_2,2,L_2)$, where $J_1,J_2,L_1,L_2$ are disjoint sequences. Then
$$\sigma\tau=(1,J_2,L_1)(2,L_2,J_1),$$
which is a cycle if and only if $J_2$ and $L_1$ are both empty or $J_1$ and $L_2$ are both empty.
We have $\sigma=(1,\ldots,2),\tau=(1,2,\ldots)$ in  the former case and $\sigma=(1,2,\ldots),\tau=(1,\ldots,2)$ in the latter case.

\qed

Denote by  $G(\mathcal{C}_n(k))$ the subgraph of $G(\Omega_n)$ induced by
$\mathcal{C}_n(k)$. Then a subset of $\Sym(n)$ is a maximum $k$-cycle clique if and
only if it is a maximum clique of $G(\mathcal{C}_n(k))$.

Denote by $ {[n]\choose r}$ the set of all $r$-subsets of $[n]$. If every pair of elements in $\mathcal{F}\subseteq {[n]\choose r}$ have nonempty intersection, then we say $\mathcal{F} $ is {\it intersecting}. We will need the following result  from   extremal set theory; see \cite[Theorem 3]{HM}.

\begin{lemma}[\cite{HM}]\label{le42}
If $2\le r\le n/2$ and $\mathcal{F}\subseteq {[n]\choose r}$. If $\mathcal{F}$ is intersecting and
$$ \lvert\mathcal{F}\rvert>{n-1\choose r-1}-{n-r-1\choose r-1}+1,$$
then $\bigcap  \mathcal{F}\ne \emptyset$.
\end{lemma}

\begin{theorem}\label{pr41}
Let $n\ge 5$ be an integer. Then the set
$$T_1(n)=\{(k ,n):1\le k\le n-1\}$$
is the unique maximum clique  in $G(\mathcal{C}_n(2))$ up to permutation similarity.
\end{theorem}

\noindent\it Proof. \rm
 For  any distinct $k_1,k_2\in\{1,2,\ldots,n-1\}$, $(k_1 ,n)$ and $(k_2 ,n)$ are adjacent in $G(\Omega_n)$, since $(k_1 ,n)^{-1}(k_2 ,n)=(k_2,k_1, n)$.  Therefore, $T_1(n)$ is a  clique  of $G(\mathcal{C}_n(2))$.

Now suppose $K$ is a maximum clique  in $G(\mathcal{C}_n(2))$, which has a cardinality of at least $n-1$. Then applying Lemma \ref{le42}, we have $\bigcap K\ne \emptyset$. Since every pair of permutations in $K$ has exactly one common number, $K$ is isomorphic to $T_1(n)$.
\qed

If a graph $G$ does not contain a triangle as a subgraph, then $G$ is called {\it triangle-free}. Denote by $K_{s,t}$ the complete bipartite graph whose vertex parts have cardinalities $s$ and $t$. The  following result is well known in   extremal graph theory.
\begin{lemma}[\cite{WM}]\label{le41}
Let $G$ be a triangle-free graph on $n$ vertices. Then
$$e(G)\le \lfloor n^2/4\rfloor$$
with equality if and only if $G$ is isomorphic to $K_{\lfloor n/2\rfloor,\lceil n/2\rceil}$.
\end{lemma}

Let $K$ be a subset of $\Sym(n)$. We denote by $K^{-1}=\{\alpha^{-1}: \alpha\in K\}$.

\begin{lemma}\label{le43}
Suppose $\mathcal{A}\subseteq\{(i,j,n):1\le i,j\le n-1\}$ is a clique of $G(\mathcal{C}_n(3))$ which  contains no pair of $(i,j,n)$ and $(j,i,n)$.
Then
\begin{equation}\label{eq41}
\mid \mathcal{A}\mid \le \lfloor (n-1)^2/4\rfloor.
\end{equation}
Moreover, if $n$ is odd,  equality in (\ref{eq41}) holds if and only if $\mathcal{A}$ is isomorphic to
$$T_2(n)=\{(i,j ,n): 1\le i \le \lfloor n/2\rfloor, \lfloor n/2\rfloor+1\le j\le n-1\};$$
if $n$ is even,  equality in (\ref{eq41}) holds if and only if $\mathcal{A}$ or $\mathcal{A}^{-1}$ is isomorphic to $T_2(n)$.
\end{lemma}

\noindent\it Proof. \rm
Construct a graph $G$ with vertex set $\{1,2,\ldots,n-1\}$ and edge set $$E(G)=\{(i,j): (i,j,n)\in \mathcal{A}\text{ or }(j,i,n)\in \mathcal{A}\}. $$
We conclude that $G$ is triangle-free. Otherwise, suppose $G$ contains a triangle $ijk$, i.e., $(i,j),(j,k),(k,i)\in E(G)$. Then $\mathcal{A}$ contains three permutations $\sigma_1,\sigma_2,\sigma_3$ such that
 $\sigma_1\in \{(i,j,n),(j,i,n)\},\sigma_2\in \{(j,k,n),(k,j,n)\},\sigma_3\in \{(i,k,n),(k,i,n)\}$. It follows that any choice of $\sigma_1,\sigma_2,\sigma_3$  generates two nonadjacent permutations in $\mathcal{A}$. Hence, $G$ is triangle-free.

 Applying Lemma \ref{le41}, we have $e(G)\le  \lfloor (n-1)^2/4\rfloor$ with equality if and only if $G$ is isomorphic to  $K_{\lfloor (n-1)/2\rfloor,\lceil (n-1)/2\rceil}$.
 Therefore, we have (\ref{eq41}). Moreover,  equality in (\ref{eq41}) holds if and only if
 $$ \mathcal{A}=\{(i,j,n): i\in S, j\in [n-1]\setminus S\} $$
with $S\subseteq [n-1]$ and $|S|=\lfloor (n-1)/2\rfloor$ or $\lceil(n-1)/2\rceil$. Therefore, if $n$ is odd, equality in (\ref{eq41}) holds if and only if
$ \mathcal{A}$ is isomorphic to $T_2(n)$. If $n$ is even,  equality in (\ref{eq41}) holds if and only if $\mathcal{A}$ or $\mathcal{A}^{-1}$ is isomorphic to $T_2(n)$.
\qed

\begin{theorem}\label{th27}
Let $n\ge 10$ be an integer. Then
$$\omega(G(\mathcal{C}_n(3)))=\lfloor (n-1)^2/4\rfloor.$$
Moreover, if $n\ge 11$ is odd, then  $T_2(n)$ is the unique maximum clique in   $G(\mathcal{C}_n(3))$ up to permutation similarity; if $n\ge 12$ is even, then $K$ is a maximum clique in $G(\mathcal{C}_n(3))$ if and only if   $K$ or $K^{-1}$ is isomorphic to $T_2(n)$.
\end{theorem}

\noindent\it Proof. \rm For $n = 10,11,12$, the results can be done by a computer search using \GAP \cite{GAP} and \Cliquer \cite{NO}; see Appendix \ref{app:max3clique}.
We assume that $n \ge 13$ in the following.
Suppose $K$ is a maximum clique in $G(\mathcal{C}_n(3))$. By Lemma \ref{le43} we have $$|K|\ge \lfloor (n-1)^2/4\rfloor.$$
Firstly we prove the following claim.

\noindent{\it {\bf Claim 1.} $K$ does not contain three distinct pairs of the form $\{\sigma,\sigma^{-1}\}.$}

Note that if $\sigma,\sigma^{-1}\in K$, then by Proposition \ref{p46}, each permutation $\delta\in K\setminus\{\sigma,\sigma^{-1}\}$ has exactly one common number with $\sigma$.

Suppose $K$ contains three distinct pairs $\{\sigma_1,\sigma_1^{-1}\},\{\sigma_2,\sigma_2^{-1}\},\{\sigma_3,\sigma_3^{-1}\}.$
Denote by
$$R=K\setminus \{\sigma_1,\sigma_2,\sigma_3,\sigma_1^{-1},\sigma_2^{-1},\sigma_3^{-1}\}.$$Without loss of generality, we may assume $\sigma_1=(1,2,3)$, $\sigma_2=(1,4,5)$. Then, there are two possibilities for $\sigma_3$:

{\it Case 1.} $\sigma_1,\sigma_2,\sigma_3$ do not have a common number, i.e., $\sigma_3$ does not contain the number 1. Then without loss of generality, we may assume $\sigma_3=(2,4,6)$.
Given any $(i,j,k)\in R$, we have
\begin{eqnarray}
&& |\{i,j,k\}\cap \{1,2,3\}|=1,\label{eq42}\\
&& |\{i,j,k\}\cap \{1,4,5\}|=1,\label{eq43}\\
&& |\{i,j,k\}\cap \{2,4,6\}|=1.
\end{eqnarray}
Hence, the set $\{i,j,k\}$ contains exactly one  of $\{1,6\}, \{2,5\}$, $\{3,4\}$  and $\{3,5,6\}$ as a subset. Therefore,
$$R \subseteq \{(3,5,6),(3,6,5),(1,6,k),(6,1,k),(2,5,k),(5,2,k),(3,4,k),(4,3,k), 7\le k\le n\} .$$
Moreover, if  $R$  contains a pair of $(1,6,k)$ and $(6,1,k)$, then it contains no other permutation with the form $(1,6,i)$ or $(6,1,i)$. Similar properties hold for $(2,5,k),(5,2,k)$ and $(3,4,k),(4,3,k)$.
It follows that
$$\lvert K\rvert = 6 + |R| \le 6+3(n-6)+2<\lfloor (n-1)^2/4\rfloor,$$
a contradiction.

{\it Case 2.} $\sigma_1,\sigma_2,\sigma_3$  have a common number. If $K$ contains another pair $\sigma_4,\sigma_4^{-1}$ which does not contain the number 1, then by Case 1 we can deduce a contradiction. So we may assume all pairs of $\{\sigma,\sigma^{-1}\}$ in $K$ have a common number 1. Suppose $K$ contains $r$ pairs of $\{\sigma,\sigma^{-1}\}$, where $r \ge 3$. Recall that any two permutations from these pairs have exactly one common number if they are not from the same pair. We may assume these pairs to be
$\{\sigma_i,\sigma_i^{-1}\}=\{(1,2i,2i+1),(1,2i+1,2i)\},i=1,2,\ldots,r$.

Given any $(i,j,k)\in R$, we have (\ref{eq42}), (\ref{eq43}) and
\begin{eqnarray}
&&\mid \{i,j,k\}\cap \{1,6,7\}\mid=1. \label{eq45}
\end{eqnarray}
If $1\not\in  \{i,j,k\}$, then by  (\ref{eq42}), (\ref{eq43}) and  (\ref{eq45}) we have $ \mid\{i,j,k\}\cap \{2,3\}\mid=1$, $ \mid\{i,j,k\}\cap \{4,5\}\mid=1$   and $ \mid\{i,j,k\}\cap \{6,7\}\mid=1$.
If $1 \in  \{i,j,k\}$, then $(i,j,k)$ has the form $(1,u,v)$ with $8\le u,v\le n$. Partition the set $K$ as $K=R_1\cup R_2\cup R_3$, where $R_1=\{\sigma_1,\sigma_1^{-1},\ldots,\sigma_r,\sigma_r^{-1}\},$ every permutation in $R_2$ contains the number 1, while every permutation in $R_3$ does not. Then $|R_3|\le 8$. Notice that $R_2$ contains no pair of $\{\sigma,\sigma^{-1}\}$, and all the permutations in $R_2$ have a common number 1. Moreover, all numbers of the permutations in $R_2$ are from $[n]\setminus\{2,3,\ldots,2r+1\}$. So $R_2$ is isomorphic to a subset of $G(\mathcal{C}_3(n-2r))$ which satisfies the condition of Lemma \ref{le43}. Hence,
$|R_2|\le \lfloor (n-2r-1)^2/4\rfloor$. Therefore,
$$\lvert K\rvert \le |R_1|+|R_2|+|R_3|\le 2r+\lfloor (n-2r-1)^2/4\rfloor+8<\lfloor (n-1)^2/4\rfloor,$$
a contradiction. This completes the proof of Claim 1.

Let $$T=\big\{\{i,j,k\}: (i,j,k)\in K\big\}.$$
By Claim 1 we have
$$|T|\ge |K|-2>{n-1\choose 2}-{n-4\choose 2}+1$$
for $n\ge 13$. Applying Lemma \ref{le42} we have $\bigcap T\ne\emptyset$.

If $\mid \cap T\mid =2$, say, $\cap T=\{1,2\}$,  then by Proposition \ref{p46}, every permutation in $K$ has the form $(1,2, k)$ or $(2,1,k)$. It follows that $|K|\le 2(n-2)<\lfloor (n-1)^2/4\rfloor$, a contradiction.
Therefore, we have $\mid \cap T\mid =1$, i.e., all permutations in $K$ have a common number, say, $n$. Hence every permutation in $K$ has the form $(i,j,n)$.

Now we claim that $K$ contains no pair of $(i,j,n)$ and $( j,i,n)$. Otherwise suppose $K$ has $k$ of  such pairs. By Claim 1 we know $k\le 2$.  If $k=2$,  without loss of generality we may assume $W\equiv\{(1,2, n), (2,1, n),(3,4, n),(4,3, n)\}\subseteq K$.
Then every permutation $\sigma\in K\setminus W$ has the form $( u,v,n)$ with $u,v\in \{5,6,\ldots,n-1\}$. So $K\setminus W$ is isomorphic to a subset of $G(\mathcal{C}_3(n-4))$ which satisfies the condition of Lemma \ref{le43}. Therefore,
$|K\setminus W|\le \lfloor(n-5)^2/4\rfloor$. It follows that
$$|K|=|W|+|K\setminus W|\le4+\lfloor(n-5)^2/4\rfloor<\lfloor(n-1)^2/4\rfloor.$$
 If $k=1$, applying similar arguments as above we have
$$|K|\le2+\lfloor(n-3)^2/4\rfloor<\lfloor(n-1)^2/4\rfloor.$$ In both cases we get  contradictions.

Finally, applying Lemma \ref{le43} on $K$ we  have
$
\lvert K \rvert \le  \lfloor (n-1)^2/4\rfloor,
$ and hence $
\lvert K \rvert= \lfloor (n-1)^2/4\rfloor.
$
Moreover, if $n\ge 11$ is odd, then $K$  is isomorphic to $T_2(n)$; if $n\ge 12$ is even, then $K$ or $K^{-1}$ is isomorphic to $T_2(n)$.
This completes the proof.
\qed

\medskip\noindent
\begin{remark} By a numerical computation, we can determine the maximum 3-cycle
cliques for $n = 3, \dots, 12$; see Appendix \ref{app:max3clique}.  For $n \le 9$, the size of maximum 3-cycle cliques is
larger than $\lfloor (n-1)^2/4\rfloor$.
\begin{center}
\begin{tabular}{|c|c|c| c| c| c| c| c|}
\hline 
$n$ & 3 & 4 & 5 & 6 & 7 & 8 & 9\cr
\hline
$\lfloor (n-1)^2/4\rfloor$ & 1 & 2 & 4 & 6 & 9 & 12 & 16\cr\hline
$\omega(G(\mathcal{C}_n(3)))$ & 1 & 2& 5&8 & 14&14 &17
\cr\hline
\end{tabular}
\end{center}
For discussions on maximum $k$-cycle cliques for $k\ge 4$, see Appendix \ref{app:maxn567} and Appendix \ref{app:addkcycle}.
\end{remark}

\begin{remark}
The celebrated Erd\H{o}s-Ko-Rado Theorem \cite{EKR} in extremal set theory states that if $K$ is a family of distinct $r$-element subsets of $[n]$ with $n\ge 2r$ such that every pair of subsets has nonempty intersection, then $|K|\le {n-1\choose r-1}$, with equality if and only if all subsets in $K$ contain a fixed common element. Similarly,  Theorem \ref{th27} serves as an analogue of the
Erd\H{o}s-Ko-Rado  Theorem for $\Sym(n)$. Specifically, it states that if $K$ is a subset of $\Sym(n)$ consisting of 3-cycle permutations such that $\delta_1^{-1}\delta_2$ is a single cycle for all $\delta_1,\delta_2\in K$, then   $K\le \lfloor (n-1)^2/4\rfloor$.
\end{remark}

\medskip
The following result provides a lower bound on $\omega(n)$.
\begin{theorem}\label{pr42}
Let $n\ge 10$ be an integer. Then
 $\{I_n\}\cup T_1(n)\cup T_2(n)$ is a maximal clique in  $G(\Omega_n)$.

\noindent
Consequently, $\omega(n)\ge (n^2+2n)/4$ when $n$ is even and $\omega(n)\ge (n+1)^2/4$ when $n$ is odd.
\end{theorem}
\noindent
\it Proof. \rm
 ~For any $(n,k)\in T_1(n)$ and $(i,j, n)\in T_2(n)$, we have
$$(k, n)(i,j, n)=\left\{\begin{array}{ll}
(i,j,k, n),&\text{if } k\not\in\{i,j\};\\
(i,j),&\text{if } k=i;\\
(i, n),& \text{if } k=j.
\end{array}\right.$$
Therefore,  $\{I_n\}\cup T_1(n)\cup T_2(n)$ is a   clique in  $G(\Omega_n)$. Now it suffices to prove that for any cycle $(i_1,i_2,\ldots,i_m)\in \Sym(n)\setminus (\{I_n\}\cup T_1(n)\cup T_2(n))$, there exists $\sigma\in   T_1(n)\cup T_2(n)$ such that
$\sigma^{-1}(i_1,i_2,\ldots,i_m)$ is not a cycle.

By Theorem \ref{pr41} and Theorem \ref{th27}, we may assume $m\ge 4$.  If $n\in \{i_1,\ldots,i_m\}$, without loss of generality, we may assume $i_m=n$. Then $$(i_2,n)(i_1,\ldots,i_{m-1},n)=(i_1,n)(i_2,\ldots,i_{m-1})$$ is not a cycle.

 Suppose $n\not\in \{i_1,\ldots,i_m\}$. If $m<n-1$, then for every $k\in [n-1]\setminus \{i_1,\ldots,i_m\}$, the permutation $(k,n)(i_1,i_2,\ldots,i_m)$ is not a cycle. If $m=n-1$, without loss of generality we may assume $i_1=1$. Notice that there exists some $t$ such that $2<t<n-1$ and $i_t\ge \lfloor n/2\rfloor+1$. We have $(1,i_t,n)\in T_2(n)$ and
$$(1,i_t,n)^{-1}(1,i_2,\ldots,i_{n-1})=(1,i_2,\ldots,i_{t-1})(i_t,\ldots,i_{n-1},n)$$
is not a cycle.

Therefore, $\{I_n\}\cup T_1(n)\cup T_2(n)$ is a maximal clique in  $G(\Omega_n)$.
\qed

For the upper bound on $\omega(n)$, we have the following theorem.

\begin{theorem}
Suppose $\mathcal{I}$ is a subgroup of $\Sym(n)$ such that it is an independent set of  $G(\Omega_n)$. Then we have
$$\omega(n)\le n!/|\mathcal{I}|.$$
In particular, if $n=2k$ be an even number, then
$$\omega(n)\le (2k)!/(k!2^{k-1}).$$

\end{theorem}
\noindent\it Proof. \rm
Since $\mathcal{I}$ is a subgroup of $\Sym(n)$, $\Sym(n)$ can be partitioned into $n!/|\mathcal{I}|$ cosets of $\mathcal{I}$. Suppose $R$ is a maximum clique of $G(\Omega_n)$. Then $R$ contains at most one permutation  from each of these cosets.
Hence, $\omega(n)\le n!/|\mathcal{I}|.$

The last assertion follows from  Theorem \ref{subgroup}.
\qed


\section{Independent subgroups}
\label{sect:subgroups}

In this section, we will construct subgroups of $\Sym(n)$ that are maximal independent sets.
We will use both abstract algebraic notation and linear algebra (matrix) notation,
and regard the elements in $\Sym(n)$ as permutation matrices when the meaning is clear.

Given a 0-1 matrix $A$, its digraph $D(A)$ is the digraph with vertex set $\{1,2,\ldots,n\}$ and arc set $\{(i,j): a_{ij}=1, 1\le i,j\le n\}$. For a permutation matrix $P$, it is clear that $D(P)$ is the disjoint union of (directed) cycles and loops. If $D(P)$ has exactly one   cycle with length larger than one, we say $P$ is a {\it cycle (matrix)}.

For any permutation $\sigma\in \Sym(n)$, we define a permutation matrix $P$ with $P_{ij}=1$ when $j=\sigma(i)$ and $P_{ij}=0$ otherwise.
For $\sigma,\tau\in \Sym(n)$ with $P$, $Q$ their permutation matrices, $\sigma^{-1}\tau$ is a cycle if and only if $P^TQ$ is a cycle matrix.
Therefore, a permutation set $T$ is an independent set in $G(\Omega_n)$ if and only if their permutation matrices form an independent set.

Clearly, if $R$ is an independent set in $G(\Omega_n)$, then so are
$$\sigma R  = \{\sigma \tau: \tau \in R\} \qquad \hbox{ and } \qquad
R \sigma =  \{\tau \sigma : \tau \in R\}.$$
Hence, we may always assume that $I_n \in R$.
Note that $\rho \in \Sym(n)$ is a cycle if and only if $\sigma\rho\sigma^{-1}$
is a cycle for any $\sigma\in \Sym(n)$. So, $\rho_1\rho_2 \in \Sym(n)$
is a cycle if and only if $\rho_2\rho_1$ is a cycle.

For $n=3$ and $4$, we have the following proposition.
\begin{proposition}\label{h1}
We have $\alpha(3)=1$ and $\alpha(4)=4$. Moreover,
$S$  is a maximum independent set of $G(\Omega_4)$ if and only if $S$
is a coset of the subgroup
$H = \{I_4, (1,2)(3,4),(1,3)(2,4),(1,4)(2,3)\}$.
\end{proposition}

\noindent\it Proof. \rm
Since the degree of each vertex in $G(\Omega_3)$ is 5, it follows that $G(\Omega_3)=K_6$ and $\alpha(3)=1$.

Suppose $S$ is a maximum independent set of $G(\Omega_4)$.
Choose any $\delta\in S$.
Denote by $N(\delta)$ the neighbour  set of $\delta$ in  $G(\Omega_4)$, whose cardinality is $\sum_{k=2}^4 {4 \choose k} (k-1)!=20$.
Then $$ S\subseteq \Sym(4) \setminus N(\delta) =\{\delta I_4, \delta(1,2)(3,4),\delta(1,3)(2,4),\delta(1,4)(2,3)\}.$$ Since any set with the form $\{\delta I_4, \delta(1,2)(3,4),\delta(1,3)(2,4),\delta(1,4)(2,3)\}$ is an independent set. The last statement is clear.
\qed

Now we present an independent subgroup of $Sym(n)$ when $n$ is even.
We will show that the independent subgroup constructed is actually a maximal
independent set.

\begin{lemma} \label{lem-even}
 Let $\tilde{\alpha}(n)$ be the maximum size of independent sets of $G(\Omega_n)$
 in which all permutations are even. Then $\alpha(n)\le 2\tilde{\alpha}(n)$.
\end{lemma}

\noindent {\it Proof.}
Denote by $h_1(X)$ and $h_2(X)$ the number of even permutations and odd permutations in a permutation set $X$, respectively. Suppose   $S\subseteq G(\Omega_n)$ is a maximum independent set. Then  we can always assume $h_1(S)\ge h_2(S)$, since otherwise we can replace $S$ with $\sigma S$ with $\sigma$ being an arbitrary transposition $(i,j)$. Therefore,
$$\alpha(n)=h_1(S)+h_2(S)\le 2h_1(S)\le 2\tilde{\alpha}(n).$$
\vskip -.4in
\qed

\begin{theorem}  \label{subgroup}
Given an even number $n=2k \ge 4$,
$H_n = \{P \oplus P: P \in \Sym(k) \}$ and
$K_n$ be the group of matrices of the form
${\small
\begin{bmatrix} I_k-D & D \cr D & I_k-D\cr\end{bmatrix}}$
such that $D$ is a diagonal 0-1 matrix with even trace.
Then
$$G_n = H_nK_n = \{hk: h \in H_n, k \in K_n\}$$ is an independent subgroup
of the alternating group  with $|G_n| = k! 2^{k-1}$.
\end{theorem}

\it Proof. \rm
If $h \in H_k, k \in G_k$, then it is easy to show that
$hkh^{-1} = k'\in K_n$.  So, $G_n = H_nK_n
= \{hk: h \in H_n, k \in K_n\}$ is a group
and
$|G_n| = |H_n||K_n|/|H_n\cap K_n|
= k! 2^{k-1}$ as $H_n\cap K_n = \{I_n\}$.

Now, we show that  $G_n$ is an independent set.
Since $a^{-1}b \in G_n$ for any $a, b \in G_n$,
we only need to show that none of the elements in
$G_n$ is a cycle. Consider
$R = hk = {\small
\begin{bmatrix} P & 0 \cr 0 & P\cr\end{bmatrix}
\begin{bmatrix} I-D & D \cr D & I-D\cr\end{bmatrix}}$.
Suppose it is a cycle. Note that the number of fixed points
equal to $\tr(R) = 2\tr(P(I-D))$ is even.
So, the cycle must be even, which implies $\det(R)=-1$. On the other hand,
it is clear that $\det(A) = 1$ for all $A \in H_n$ or $A \in K_n$, which leads to $\det(R)=1$, a contradiction.
\qed

We will show that $G_n$ constructed in the previous theorem
is a maximal independent set in $\Sym(n)$, i.e.,
$G_n \cup \{\tau\}$ is not an independent set for any $\tau \in \Sym(n) \setminus G_n$. We will prove that every nontrivial  (right)
coset of $G_n$ contains a cycle.  First, we will need some preliminary results and notation.

For $j \in \{1, \dots, n\}$, $n = 2k$, we define
\[j^c := \begin{cases}
          j + k, \text{ if } j \le k,\\
          j - k, \text{ if } j > k.
         \end{cases}
 \]
Using this notation, $G_n$ stabilizes the partition $\cP$ of $\{1,\dots, n\}$
into the sets $\{j, j^c\}$.  That is, $\sigma$
will induce a permutation on $\cP$. Clearly, the stabilizer group $X_n$ of
$\cP$ consists of $\phi \in \Sym(n)$
that permute the subsets in $\cP$ and maps $\{j,j^c\}$ to itself for each $j$. So, it is the wreath
product of $\Sym(n/2)$ and $\Sym(2)$, and the group $G_n$ is just $X_n \cap \Alt(n)$.

Define \[G_{n,j,j^c} := G_n \cap {\rm Sym}(\{1,\dots, n\} \setminus \{j,j^c\}).\]
In other words, $G_{n,j,j^c}$ stabilizes $\mathcal{P}$, fixes both $j$ and $j^c$,
and is isomorphic to $G_{n-2}$.  We will need one lemma.

\begin{lemma}
 \label{lem:i->j}
 Assume that every nontrivial coset of $G_{n-2}$ in $\Sym({n-2})$ contains a cycle.
 If $\sigma \in \Sym(n) \setminus G_n$ and $\{i, i^c\}^\sigma = \{j, j^c\}$ for some (not necessarily distinct)
 $i,j \in \{1,\dots,n\}$, then the coset
 $\sigma G_n$ contains a cycle $\rho$ that fixes both $j$ and $j^c$.
\end{lemma}

\noindent \it Proof. \rm
 We will prove this in a few steps.  First, we have
 $$\sigma \notin G_{n,j,j^c} \le {\rm Sym}(\{1,\dots, n\} \setminus \{j,j^c\}).$$
 Moreover,
 \[(j , n/2)(j^c , n)G_{n,j,j^c}((j , n/2)(j^c , n))^{-1} = G_{n-2}\] and
 \[(j , n/2)(j^c , n){\rm Sym}(\{1,\dots, n\} \setminus \{j,j^c\})((j , n/2)(j^c , n))^{-1} = \Sym({n-2}).\]
 Since every nontrivial coset of $G_{n-2}$ in $\Sym({n-2})$ contains a cycle,
 the coset $\sigma G_{n,j,j^c} \subseteq \sigma G_n$ contains a cycle that fixes both $j$ and $j^c$.

 Next, assume $i = j$.  If $\sigma(j) = j$, $\sigma(j^c) = j^c$, we are done by the argument above.  Otherwise, let $\ell \not\in \{j, j^c\}$.  Then, $\gamma = (\ell , \ell^c)(j , j^c) \in G_n$, and $\sigma \gamma(j) = j$, $\sigma \gamma(j^c) = j^c$.  By the argument in the previous paragraph, $\sigma \gamma G_n = \sigma G_n$ contains a cycle fixing both $j$ and $j^c$.

 Finally, assume $i \neq j$.  Let $\gamma = (i , j)(i^c , j^c) \in G_n$.
 Then, $\sigma\gamma(\{j, j^c\})= \{j, j^c\}$. By the argument in the previous paragraph, $\sigma \gamma G_n = \sigma G_n$ contains a cycle fixing both $j$ and $j^c$, as desired.
\qed

Now,  we are ready to prove that $G_n$ is a maximal independent set.

\begin{theorem}
 \label{thm:Gnind}
 Let $n \ge 4$ be an even integer.  Then  $G_n$ is a maximal
 independent set in $\Sym(n)$.
\end{theorem}

\noindent \it Proof. \rm
 We know $G_n$ is an independent set by Theorem \ref{subgroup}.  We will prove that every nontrivial coset of $G_n$ contains a cycle.

Let $k=n/2$. We proceed by induction on $k$.  For $n = 4$, the result can be verified easily by Proposition \ref{h1}.
We assume now that the result is true for $k-1$ (and, hence, $n-2$).  Let $\sigma \in \Sym(n) \setminus G_n$.
We may assume for all $i \in \{1, \dots, n\}$ that $\sigma(i) = j$ implies $\sigma(i^c) \neq j^c$; otherwise, we are done by Lemma \ref{lem:i->j}.

 Suppose $\sigma(k) = i$ and $\sigma(k^c) = \sigma(n) = j^c$.  Again, if $i = j$, we are done.
 Let $\gamma = (i , j)$.  Then, $\gamma \sigma (k) = j$ and $\gamma \sigma(n) = j^c$.
 By Lemma \ref{lem:i->j}, $\gamma \sigma G_n$ contains a cycle $\rho = (a_1 , a_2, \dots , a_m)$
 fixing both $j$ and $j^c$ (that is, $j, j^c \notin \supp(\rho)$).  Moreover,
 $\sigma G_n = \gamma^{-1} \rho G_n = \gamma \rho G_n$.

 We now divide it into two cases.  Suppose first $i \in \supp(\rho)=\{a_1, \dots, a_m\}$.
 Since $j \notin \supp(\rho)$, $\gamma \rho$ itself is a cycle, and we are done.

 Otherwise, suppose $i \notin \supp(\rho)$.  In this case, since $\supp(\rho) \cap \{i,j\}
 = \varnothing$, we have $\rho \gamma = \gamma \rho$.  If $i^c \in \supp(\rho)$, then,
 letting $\beta = (i , j)(i^c , j^c) \in G_n$, we have $\rho \gamma \beta = \rho (i^c , j^c)$,
 $i^c \in \supp(\rho)$, $j^c \notin \supp(\rho)$, and we are done as above.

 So, we may assume $\supp(\rho) \cap \{i,i^c, j, j^c\} = \varnothing$.
If $a_m^c \in \supp(\rho)$, i.e., if $a_m^c = a_\ell$ for some $\ell$, $1 \le \ell < m$, then,
if $\beta = (i^c , j^c , a_\ell)(i , j , a_m) \in G_n$, we have $\rho \gamma \beta
\in \rho \gamma G_n = \sigma G_n$ and
 \[ \rho \gamma \beta = (a_1 , \dots , a_\ell , i^c , j^c , a_{\ell + 1} \dots , a_m , j).\]
 Otherwise, $a_m^c \notin \supp(\rho)$, in which case, if $\beta = (j , j^c)(i , a_m , i^c , a_m^c)=(i,i^c)(j,j^c)(i,a_m)(i^c,a_m^c) \in G_n$ (since $\beta$ preserves $\mathcal{P}$ and is an even permutation), $\rho \gamma \beta \in \sigma G_n$, and
 \[ \rho \gamma \beta = (a_1 , \dots , a_m , i^c , a_m^c , j , j^c , i).\]
 Therefore, in any case, $\sigma G_n$ always contains a cycle when $\sigma \notin G_n$, and $G_n$ is therefore a maximal independent set.
\qed

In fact, even more is true:

\begin{theorem}
 \label{thm:GnSn+1}
 Let $n  \ge 4$ be an even integer.  If we identify the group
 $G_n$ as a subgroup of $\Sym(n+1)$ with elements fixing $n+1$,
 then $G_n$  is a maximal independent set in $\Sym({n+1})$.
\end{theorem}

\noindent \it Proof. \rm
 Let $\sigma \in \Sym({n+1})$.
 If $\sigma(n+1) = n+1$, then $\sigma \in \Sym(n)$, and the result holds by Theorem \ref{thm:Gnind}.  Hence we may assume $\sigma(n+1) = j < n+1$.  If $\gamma = (j , n+1)$, then $\gamma \sigma(n+1) = n+1$, and, by Theorem \ref{thm:Gnind}, $\gamma \sigma G_n = \rho G_n$, where $\rho \in \Sym(n)$ (i.e., $n+1 \notin \supp(\rho)$) is a cycle.  Thus, $\sigma G_n = \gamma \rho G_n$.  For $i \in \{1,\dots, n\}$, we will use the notation $i^c$ for the element in $\{1,\dots,n\}$ such that $\{i,i^c\}$ is a block preserved by $G_n$.

 If $j \in \supp(\rho)$, then $\gamma \rho$ is a cycle, and we are done.  Otherwise, $\supp(\rho) \cap \{j, n+1\} = \varnothing$, and so $\rho \gamma = \gamma \rho$.  Let $\rho = (a_1 , a_2 , \dots , a_m)$.  Assume first $j^c \in \supp(\rho)$.  Up to relabeling, we may assume $j^c = a_m$.  If $a_1^c \notin \supp(\rho)$, then, if $\beta = (a_1 , a_1^c)(j , j^c) \in G_n$, we have $\rho \gamma \beta \in \gamma \rho G_n = \sigma G_n$ and
 \[ \rho \gamma \beta = (a_1 , a_1^c , a_2 , \dots , a_{m-1} , j^c , n+1 , j), \]
 and we are done.  If $a_1^c = a_r \in \supp(\rho)$, then, if $\beta =(a_1 , j^c)(a_1^c , j)   \in G_n$, then $\rho \gamma \beta \in \sigma G_n$,
 \[\rho \gamma \beta = (a_2 , \dots , a_{r-1} , a_1^c , n+1 , j , a_{r+1} , \dots , a_{m-1} , j^c),\]
 and we are done.

 So, we may assume $j, j^c \notin \supp(\rho)$.  If $a_1^c \in \supp(\rho)$, say $a_1^c = a_r$, then, for $\beta = (a_1 , j^c)(a_1^c , j) \in G_n$, we have $\rho \gamma \beta \in \sigma G_n$,
 \[ \sigma \gamma \beta = (a_1 , j^c , a_2 , \dots , a_{r-1} , a_1^c , n+1 , j , a_{r+1} , \dots, a_m),\]
 and we are done.  Since we may reorder the indices, we can thus assume $a_i^c \notin \supp(\rho)$ for all $a_i \in \supp(\rho)$; in particular, $a_1^c, a_m^c \notin \supp(\rho)$.  Thus, if $\beta = (a_1^c , a_m , j^c)(a_1 , a_m^c , j) \in G_n$, then $\rho \gamma \beta \in \sigma G_n$,
 \[ \rho \gamma \beta = (a_1 , a_m^c , n+1 , j , a_2 , \dots , a_m , j^c , a_1^c).\]
 This exhausts all the cases, so every nontrivial coset of $G_n$ in $\Sym({n+1})$ always contains a cycle, and, therefore, $G_n$ is a maximal independent set in $\Sym({n+1})$.
\qed

We remark that $G_n$ is not a maximum independent set in
$\Sym({n+1})$  in general.

\section{Related topics and further research}
\label{sect:related}

The transportation polytope  is defined  as
\begin{eqnarray*}P_{m,n}(r,s)=\left\{\begin{array}{ll}
 &\sum_{j=1}^nx_{ij}=r_i, i=1,2,\ldots,m; \\
X=(x_{ij})_{m\times n}: &\sum_{i=1}^mx_{ij}=c_j,j=1,2,\ldots,n;\\
&x_{ij}\ge 0, i=1,2,\ldots,m,j=1,2,\ldots,n\end{array}\right\},
\end{eqnarray*}
where $r=(r_1,r_2,\ldots,r_m)$ and $s=(c_1,c_2,\ldots,c_n)$ are positive vectors. The transportation polytope arises from the transportation problem, which was posed by Hitchcock \cite{FH} and Koopmans \cite{TK} and is well known in   operation research and mathematical programming. The transportation polytope attracts lots of research; see \cite{AB,EB,BVS,DKOS,FL,JO} and references therein.

When $m=n$ and $r=s=(1,1,\ldots,1)$, the transportation polytope $P_{m,n}(r,s)$,
reduces to the Birkhoff Polytope $\Omega_n$. It is interesting to see whether our
study can be extended to the transportation polytope.

Also, one may consider the extension of doubly stochastic matrices to doubly stochastic
tensors, i.e.,  $(a_{ijk})$ such that $a_{ijk} \ge 0$ with
$1 \le i,j,k\le n$, and $\sum_{i} a_{ijk} = 1$ for any fixed $(j,k)$ pairs,
$ \sum_{j} a_{ijk} = 1$ for any fixed $(i,k)$ pairs, and
$\sum_{k} a_{ijk} = 1$ for any fixed $(i,j)$ pairs.
It would be interesting to determine the vertices of the convex set of
doubly stochastic tensors, and define a suitable Cayley graph on the vertices.

\section*{Acknowledgment}

The research of Huang was supported by the National Natural Science Foundation of China
(No. 12171323).
Li is an affiliate member of the Institute for Quantum Computing, University
of Waterloo; his
research was partially supported by the Simons Foundation Grant 851334.
The research of Sze was supported by a Hong Kong RGC grant PolyU 15305719 and a PolyU Internal Research Fund 4-ZZKU.

\newpage
\appendix
\section{Appendix}
\label{appendix}

The programs \GAP \cite{GAP} and \Cliquer \cite{NO} were used extensively for the calculations done throughout Appendix \ref{appendix}.

\subsection{Maximum cliques attaining $\omega(n)$ for $n=5,6,7$}
\label{app:max567}

For $n=5$, the following set is a maximum clique in $G(\Omega_5)$:
 $\{I_5$,           (1,2,3,4,5),  (1,2,4,5,3),  (1,2,5,3,4),
(1,3,5,4,2), (1,3,2,5,4),  (1,3,4,2,5),  (1,4,3,5,2),
(1,4,5,2,3), (1,4,2,3,5),  (1,5,4,3,2),  (1,5,2,4,3), (1,5,3,2,4)$\}$.

For $n=6$, the following set is a maximum clique in $G(\Omega_6)$:
$\{I_6$,     (3,4),  (3,5,4),   (3,6,4),  (2,4,3), (1,4,3),
(1,2,3,4,6,5), (1,2,5,3,4,6), (1,2,6,5,3,4), (1,3,4,5,6,2),
(1,3,4,2,5,6), (1,3,4,6,2,5), (1,5,6,3,4,2), (1,5,2,6,3,4),
(1,5,3,4,2,6), (1,6,3,4,5,2), (1,6,5,2,3,4), (1,6,2,3,4,5)$\}$.

For $n=7$, the following set is a maximum clique in $G(\Omega_7)$:
$\{I_7$, (4,5), (2,4,3), (2,4), (2,4,6), (2,4,7), (2,5,3), (2,5,6), (2,5,7), (1,2,5), (1,3,6,7,5,4,2),
      (1,3,5,4,2,6,7), (1,3,7,5,4,2,6), (1,4,5), (1,5,4,2,3,7,6), (1,5,4,2,6,3,7), (1,5,4,2,7,6,3),
      (1,6,7,3,5,4,2), (1,6,5,4,2,7,3), (1,6,3,5,4,2,7), (1,7,3,6,5,4,2), (1,7,6,5,4,2,3), (1,7,5,4,2,3,6)$\}$.\\

\subsection{Maximum $3$-cycle cliques for $n \le 12$}
\label{app:max3clique}

For $n = 3$, it is clear that a maximum $3$-cycle clique consists of a single $3$-cycle $(1,2,3)$, and this is unique up to permutation similarity.

For $n = 4$, up to permutation similarity there are two maximum $3$-cycle cliques of size $2$, namely $\{(2,3,4), (1,4,2) \}$ and $\{ (2,3,4), (2,4,3)\}$.

For $n = 5$, there are two maximum $3$-cycle cliques of size $5$ up to permutation similarity:
\begin{eqnarray*}
&& \{ (2,5,3), (2,5,4), (1,2,5), (1,3,4), (1,4,3)\},\\
&& \{ (3,4,5), (2,3,4), (1,2,3), (1,2,5), (1,4,5)\}.
\end{eqnarray*}

For $n = 6$, there are six maximum $3$-cycle cliques of size $9$ up to permutation similarity:
\begin{eqnarray*}
&& \{(3,5,6), (2,4,5), (2,5,4), (2,6,3), (1,3,2), (1,3,5), (1,4,6), (1,6,4)\}, \\
&& \{ (4,5,6), (4,6,5), (2,3,4), (2,3,5), (2,3,6), (1,3,4), (1,3,5), (1,3,6)\}, \\
&& \{(4,5,6), (4,6,5), (2,3,4), (2,3,5), (2,3,6), (1,3,4), (1,3,5), (1,6,2)\}, \\
&& \{(4,5,6), (4,6,5), (2,3,4), (2,3,5), (2,3,6), (1,3,4), (1,5,2), (1,6,2)\} \\
&& \{(4,5,6), (4,6,5), (2,3,4), (2,3,5), (2,3,6), (1,4,2), (1,5,2), (1,6,2)\}, \\
&& \{(4,5,6), (4,6,5), (2,3,4), (2,4,3), (1,2,5), (1,3,6), (1,5,2), (1,6,3)\} .
\end{eqnarray*}

For $n = 7$, there is a unique maximum $3$-cycle clique of size $14$ up to permutation similarity:
$$\begin{array}{lllllll}
 \{ (4,6,7), &(4,7,6), &(3,4,5), &(3,5,4), &(2,3,7), &(2,5,6), &(2,6,5), \\
 \ (2,7,3), &(1,2,4), &(1,3,6), &(1,4,2), &(1,5,7), &(1,6,3), &(1,7,5)\}.
\end{array}$$

For $n = 8$, there are three maximum $3$-cycle cliques of size $14$ up to permutation similarity:
$$\begin{array}{lllllll}
\{(6,7,8), &(5,6,7), &(4,6,7), &(3,4,6), &(3,5,6), &(3,8,6), &(2,3,7), \\
\ (2,4,6), &(2,5,6), &(2,7,3), &(2,8,6), &(1,6,2), &(1,6,3), &(1,6,7)\}, \\ \\
\{ (6,7,8), &(6,8,7), &(4,5,6), &(4,5,7), &(4,5,8), &(3,6,4), &(3,7,4), \\
 \ (3,8,4), &(2,6,4), &(2,7,4), &(2,8,4), &(1,6,4), &(1,7,4), &(1,8,4)\}, \\ \\
 \{(6,7,8), &(6,8,7), &(4,5,6), &(4,6,5), &(3,4,7), &(3,5,8), &(3,7,4), \\
 \ (3,8,5), &(2,3,6), &(2,4,8), &(2,5,7), &(2,6,3), &(2,7,5), &(2,8,4)\}.
\end{array}$$

For $n = 9$, there are two maximum $3$-cycle cliques of size $17$ up to permutation similarity:
$$\begin{array}{lllllllll}
  \{ (2,3,5), &(2,3,6), &(2,4,5), &(2,4,6), &(2,7,5), &(2,7,6), &(2,8,5), &(2,8,6),
      &(2,9,5),\\
  \  (2,9,6), &(1,2,3), &(1,2,4), &(1,2,7), &(1,2,8), &(1,2,9), &(1,5,6)\},
      &(1,6,5)\}, & \\ \\
\{   (7,8,9), &(7,9,8), &(5,6,7), &(5,6,8), &(5,6,9), &(4,6,7), &(4,6,8), &(4,6,9),
       &(3,6,7),\\
   \   (3,6,8), &(3,6,9), &(2,6,7), &(2,6,8), &(2,6,9), &(1,6,7), &(1,6,8)\},
      &(1,6,9)\}.&
  \end{array}$$

For $n = 10$, there are four maximum $3$-cycle cliques of size $20$ up to permutation similarity:
$$\begin{array}{lllllll}
\{ (5,6,7), &(5,8,7), &(5,9,7), &(5,10,7),  &(4,6,7), &(4,8,7), &(4,9,7), \\
 \ (4,10,7), &(3,6,7), &  (3,8,7), &(3,9,7), &(2,10,7), &(1,7,2), &(1,7,3),\\
 \ (3,10,7), &(2,6,7), &(2,8,7), &(2,9,7), &(1,7,4), &(1,7,5)\}, &  \\ \\
\{ (8,9,10), &(8,10,9), &(6,7,8), &(6,7,9), &(6,7,10), &(5,7,8), &(5,7,9), \\
 \ (5,7,10), &(4,7,8), &(4,7,9), &(4,7,10), &(3,7,8), &(3,7,9), &(3,7,10),\\
 \ (2,7,8), &(2,7,9), &(2,7,10), &(1,7,8), &  (1,7,9), &(1,7,10)\}, & \\ \\
 \{ (8,9,10), &(8,10,9), &(6,7,8), &(6,7,9), &(6,7,10), &(5,8,6), &(5,9,6),\\
\  (5,10,6), &(4,8,6), & (4,9,6), &(4,10,6), &(3,8,6), &(3,9,6), &(3,10,6),\\
\  (2,8,6), &(2,9,6), &(2,10,6), &(1,8,6), &   (1,9,6), &(1,10,6)\}, & \\ \\
\{ (8,9,10), &(7,8,9), &(6,7,8), &(6,10,8), &(5,7,8), &(5,10,8), &(4,7,8),\\
\  (4,10,8), &(3,7,8), &(3,10,8), &(2,8,3), &(2,8,4), &(2,8,5), &(2,8,6),\\
\  (2,8,9), &(1,8,3), &(1,8,4), &(1,8,5), &1,8,6), &(1,8,9)\}. &
  \end{array}$$

For $n = 11$, there is one maximum $3$-cycle clique of size $25$ up to permutation similarity:
$$\begin{array}{lllllllll}
   \{(7,11,8), &(7,11,9), &(7,11,10), &(6,8,7), &(6,9,7), &(6,10,7), &(5,8,7), &(5,9,7), &(5,10,7),\\
   (4,7,5), &(4,7,6), &(4,7,11), &(3,4,7), &(3,8,7), &(3,9,7), &(3,10,7), &(2,7,3), &(2,7,5),\\
   (2,7,6), &(2,7,11), &(1,2,7), &(1,4,7), &(1,8,7), &(1,9,7), &(1,10,7)\}. & &
  \end{array}$$

For $n = 12$, there are two maximum $3$-cycle cliques of size $30$ up to permutation similarity:
$$\begin{array}{lllllllll}
\{(4,7,5), &(4,7,6), &(4,7,8), &(4,7,9), &(4,10,5), &(4,10,6), &(4,10,8), &(4,10,9), &(4,11,5),\\
   (4,11,6), &(4,11,8), &(4,11,9), &(4,12,5), &(4,12,6), &(4,12,8), &(4,12,9), &(3,4,7), &(3,4,10),\\
   (3,4,11), &(3,4,12), &(2,3,4), &(2,5,4), &(2,6,4), &(2,8,4), &(2,9,4), &(1,3,4), &(1,5,4),\\
   (1,6,4), &(1,8,4), &(1,9,4)\}, & & & & & &
  \end{array}$$
$$\begin{array}{llllllll}
\{ (10,11,12), &(9,10,11), &(8,9,10), &(8,12,10), &(7,9,10), &(7,12,10),
      &(6,9,10), &(6,12,10),\\
      (5,9,10), &(5,12,10), &(4,10,5), &(4,10,6), &(4,10,7), &(4,10,8), &(4,10,11), &(3,10,5),\\
      (3,10,6), &(3,10,7), &(3,10,8), &(3,10,11), &(2,10,5), &(2,10,6), &(2,10,7), &(2,10,8),\\
      (2,10,11), &(1,10,5), &(1,10,6), &(1,10,7), &(1,10,8), &(1,10,11)\}. & &
  \end{array}$$

\subsection{Maximum $n$-cycle cliques for $n=5,6,7$}
\label{app:maxn567}

For $n =  5$, the maximum $5$-cycle clique size is 12:
$$\begin{array}{llllll}
\{(1,2,3,4,5),&(1,2,4,5,3),&(1,2,5,3,4),
(1,3,2,5,4),&(1,3,4,2,5),&(1,3,5,4,2),\\
\ (1,4,2,3,5),&(1,4,5,2,3),&(1,4,3,5,2),
(1,5,2,4,3),&(1,5,3,2,4),&(1,5,4,3,2)\}.
\end{array}$$
Note that if $\sigma$ is in the clique, then $\sigma^{-1}$ is also in the clique.
The computer also shows that this clique and the identity permutation form a maximum clique of
$\Sym(5)$.

For $n =  6$, the maximum $6$-cycle clique size is also 12. The clique can be constructed from the case $n=5$ by inserting the number 6 after the number 1 in all permutations as follows.
$$\begin{array}{llllll}
\{(1,6,2,3,4,5),&(1,6,2,4,5,3),&(1,6,2,5,3,4),
(1,6,3,2,5,4),&(1,6,3,4,2,5),&(1,6,3,5,4,2),\\
\ (1,6,4,2,3,5),&(1,6,4,5,2,3),&(1,6,4,3,5,2),
(1,6,5,2,4,3),&(1,6,5,3,2,4),&(1,6,5,4,3,2)\}.
\end{array}$$
We can add $(1,6), (1,2,6), (1,3,6), (1,4,6), (1,5,6)$ and the identity permutation
to get a maximum clique with 18 elements in $\Sym(6)$.

For $n = 7$, the maximum $7$-cycle clique size is 18. The clique is unique up to
permutation similarity:
$$\begin{array}{lllll}
\{     (1,     7,     6,     5,     4,     3,     2), &
     (1,     3,     7,     6,     5,     4,     2), &
     (1,     3,     2,     7,     6,     5,     4), &
     (1,     5,     7,     6,     4,     3,     2), &
     (1,     3,     5,     7,     6,     4,     2),\\
\  (1,     3,     2,     5,     7,     6,     4), &
     (1,     2,     4,     6,     7,     5,     3), &
     (1,     4,     6,     7,     5,     2,     3), &
     (1,     2,     3,     4,     6,     7,     5), &
     (1,     2,     4,     7,     5,     6,     3),\\
\     (1,     4,     7,     5,     6,     2,     3), &
     (1,     2,     3,     4,     7,     5,     6), &
     (1,     6,     5,     7,     4,     3,     2), &
     (1,     3,     6,     5,     7,     4,     2), &
     (1,     3,     2,     6,     5,     7,     4),\\
\     (1,     2,     4,     5,     6,     7,     3), &
     (1,     4,     5,     6,     7,     2,     3), &
     (1,     2,     3,     4,     5,     6,     7)\}. &
     &
     \end{array}$$
Note that $\sigma$ in the set if and only if $\sigma^{-1}$ also in the set.\\

For $n=8$, the maximum $8$-cycle clique size is $24$, and it is unique up to permutation similarity:
$$\begin{array}{llll}
\{  (1,2,7,4,3,6,5,8), &(1,3,4,7,2,6,5,8), &(1,4,3,2,7,6,5,8), &(1,5,7,2,3,4,6,8),\\ \ (1,5,8,7,2,3,4,6), &(1,5,3,4,7,2,6,8), &(1,5,8,3,4,7,2,6), &(1,5,4,3,2,7,6,8),\\ \ (1,5,8,4,3,2,7,6), &(1,5,2,7,4,3,6,8), &(1,5,8,2,7,4,3,6), &(1,6,7,2,4,3,8,5),\\ \ (1,6,3,4,2,7,8,5), &(1,6,2,7,3,4,8,5), &(1,6,4,3,7,2,8,5), &(1,7,2,3,4,6,5,8),\\ \ (1,8,5,6,4,3,7,2), &(1,8,5,6,7,2,4,3), &(1,8,6,7,2,4,3,5), &(1,8,6,4,3,7,2,5),\\ \
  (1,8,5,6,2,7,3,4), &(1,8,6,3,4,2,7,5), &(1,8,5,6,3,4,2,7), &(1,8,6,2,7,3,4,5)\}.
  \end{array}$$

\par

\subsection
{\bf Additional observations on $k$-cycle cliques}
\label{app:addkcycle}

For general $k$, we have the following lower bound on the size of maximum $k$-cycle cliques.
\begin{proposition}
Let $n,k$ be integers such that $n\ge k\ge 4$. Then
\begin{equation}\label{eq03}
\omega(G(\mathcal{C}_n(k)))\ge\lfloor (n-k+2)^2/4\rfloor.
\end{equation}
\end{proposition}

\noindent\it Proof. \rm
 Let
 $$K(n)=\left\{(i,j ,n-k+3,n-k+4,\ldots,n): \hskip 2in \ \right.$$
 $$\left. \qquad  \qquad  1\le i \le \lfloor (n-k+2)/2\rfloor,
   \lfloor (n-k+2)/2\rfloor+1\le j\le n-k+2\right\}.$$ Note that $$(i_1,j_1 ,n-k+3, \ldots,n)(i_2,j_2 ,n-k+3, \ldots,n)^{-1}=(i_1,j_1,n-k+3)(i_2,j_2 ,n-k+3)^{-1}$$
when $\{i_1,j_1\}\bigcap\{i_2,j_2\}=\emptyset$.
Since $\sigma^{-1}\tau$ is a cycle if and only if $\sigma\tau^{-1}$ is a cycle,
we can apply   Lemma \ref{le43}
to conclude that $K$ is a clique with size
$\lfloor (n-k+2)^2/4\rfloor$.  Hence, we have (\ref{eq03}).
\qed

\medskip

\begin{remark}
For the maximum $n$-cycle cliques, when $n=4$,   the maximum clique size is 3. $\Sym(4)$ contains the following six 4-cycles:
$$
     (1,    2, 3, 4), \quad (1,4,3,2) \quad (1,2,4,3), \quad  (1,3,4,2), \quad
(1,3,2,4), \quad
      (1,4,2,3).
 $$  Note that for any 4-cycle $\sigma$, $\sigma^2$ is not a cycle. To choose the maximum clique from 4-cycles, we can pick at most one from each of the three pairs:
 $\{   (1,    2, 3, 4),  (1,4,3,2)\},$ $\{ (1,2,4,3),   (1,3,4,2)\},$ $\{(1,3,2,4),
      (1,4,2,3)\}.$ A simple calculation shows that an such arbitrary choice generates a clique. For $n=5,6,7$, by a computer search we have
$ \omega(G(\mathcal{C}_5(5)))=12, \omega(G(\mathcal{C}_6(6)))=12,
\omega(G(\mathcal{C}_7(7)))=18$ as shown in the previous subsection.
\end{remark}

\medskip\noindent
\begin{remark} If $n=7$ or 8, then
 $\{(4,5,6),(4,6,5)\}\cup T_2(7)$ and  $ \{(5,6,7),(5,7,6)\}\cup T_2(8)$ are maximal cliques in $G(\mathcal{C}_3(7)$ and $G(\mathcal{C}_3(8))$, respectively.
\end{remark}

It is interesting to note that for $n = 7$, we get a 3 cycle cliques with 15 elements:
$$ I_7, (4,6,7), (4,7,6), (3,4,5), (3,5,4), (2,3,7), (2,5,6), (2,6,5),
$$
$$\quad\ (2,7,3), (1,2,4), (1,3,6), (1,4,2), (1,5,7), (1,6,3), (1,7,5).$$
Excluding the identity, they form 7 pairs of element of the form
$(i,j,k)$ and $(i,k,j)$,
where $\{i,j,k\} \subseteq \{1, \dots, 7\}$ are the lines in the Fano plane.

Generally, for any odd prime power $q$, let $n=q^2+q+1$. Then it is known that there exists a projective plane $ \mathcal{F}=(\mathcal{P},\mathcal{L})$ of order $q$, where $\mathcal{P}$ and $\mathcal{L}$ are the sets of points and lines, respectively. Moreover, $|\mathcal{L}|=|\mathcal{P}|=n$ and each line in $ \mathcal{F}$ contains exactly $q+1$ points. For convenience, we denote the points of $ \mathcal{F}$ by $\{1,2,\ldots,n\}$.  Suppose $K$ is a $(q+1)$-cycle clique in $G(\Omega_{q+1})$.
Let $$K'=\{(i_{j_1}, \ldots,i_{j_{q+1}}): \{i_1, \ldots,i_{q+1}\} \text{  is a line in } \mathcal{F}, (j_1, \ldots,j_{q+1})\in K\}.$$
Then $K'$ is a $(q+1)$-cycle clique in $G(\Omega_{n})$ with $n|K|$ elements. In fact, for any $\sigma, \tau\in K'$, we have one of the following possibilities.
 \begin{itemize}
 \item $\sigma=(i_{j_1}, \ldots,i_{j_{q+1}}), \tau=(i_{j'_1}, \ldots,i_{j'_{q+1}})$ with $(j_1,\ldots,j_{q+1}), (j'_1,\ldots,j'_{q+1})\in K$ and $\{i_1,\ldots,i_{q+1}\}$ being a line in $ \mathcal{F}$;
 \item $\sigma=(i_{j_1},i_{j_2},\ldots,i_{j_{q+1}}), \tau=(k_{j'_1},k_{j'_2},\ldots,k_{j'_{q+1}})$ with $(j_1,\ldots,j_{q+1}), (j'_1,\ldots,j'_{q+1})\in K$, and $\{i_1,\ldots,i_{q+1}\},$ $\{k_1,\ldots,k_{q+1}\}$ being distinct lines  in $ \mathcal{F}$.
\end{itemize}
In both cases it is easy to check that $\sigma^{-1}\tau$ is a cycle.  Hence, $K'$ is a  $(q+1)$-cycle clique in $G(\Omega_n)$ with size $n|K|$. Therefore, we have
$$\omega(G(\mathcal{C}_n(q+1)))\ge n\omega(G(\mathcal{C}_{q+1}(q+1))).$$

\end{document}